\documentclass[10pt]{article}
\newcommand{\qvar}[1]{\langle {#1} \rangle}

\newcommand{\DDfDxx}[2]{\frac{d^2 #1}{d {#2}^2}}

\newcommand{\beq}{\begin{equation}}
  \newcommand{\eeq}{\end{equation}}
  \newcommand{\beqno}{\begin{displaymath}}
  \newcommand{\eeqno}{\end{displaymath}}
  \newcommand{\beqar}{\begin{eqnarray}}
  \newcommand{\eeqar}{\end{eqnarray}}
  \newcommand{\beqarno}{\begin{eqnarray*}}
  \newcommand{\eeqarno}{\end{eqnarray*}}

\newcommand{\cref}[1]{(\ref{#1})}

\newcommand{\ba}{\begin{array}}
\newcommand{\ea}{\end{array}}

\newcommand{\ben}{\begin{enumerate}}
  \newcommand{\en}{ \end{enumerate}}
\newcommand{\bei}{\begin{itemize}}
  \newcommand{\ei}{ \end{itemize}}
\newcommand{\bed}{\begin{description}}
  \newcommand{\ed}{\end{description}}
\newcommand{\bec}{\begin{center}}
  \newcommand{\ec}{\end{center}}
\newcommand{\bprop}{\begin{proposition}}
  \newcommand{\eprop}{\end{proposition}}
\newcommand{\bdf}{\begin{definition}}
  \newcommand{\edf}{\end{definition}}
\newcommand{\bth}{\begin{theorem}}
  \newcommand{\eth}{\end{theorem}}
\newcommand{\bcon}{\begin{con}}
  \newcommand{\econ}{\end{con}}
\newcommand{\bcor}{\begin{corollary}}
  \newcommand{\ecor}{\end{corollary}}
\newcommand{\bpr}{\begin{problem}}
  \newcommand{\epr}{\end{problem}}
\newcommand{\blem}{\begin{lemma}}
  \newcommand{\elem}{\end{lemma}}
 \newcommand{\bass}{\begin{assumption}}
  \newcommand{\eass}{\end{assumption}}
 \newcommand{\brem}{\begin{remark}}
  \newcommand{\erem}{\end{remark}}
\newcommand{\bres}{\begin{result}}
  \newcommand{\eres}{\end{result}}
\newcommand{\bexm}{\begin{example}}
  \newcommand{\eexm}{\end{example}}
\newcommand{\bef}{\begin{fact}}
  \newcommand{\ef}{\end{fact}}

\newtheorem{assumption}{Assumption}[section]
  \newtheorem{definition}{Definition}[section]
  \newtheorem{result}{Result}[section]
  \newtheorem{theorem}{Theorem}[section]
  \newtheorem{corollary}{Corollary}[section]
  \newtheorem{lemma}{Lemma}[section]
  \newtheorem{proposition}{Proposition}[section]
  \newtheorem{example}{Example}[section]

  \newtheorem{"definition"}{"Definition"}[section]
 \newtheorem{remark}{Remark}[section]
 \newtheorem{con}{Conjecture}[section]
 \newtheorem{problem}{Problem}[section]

\newtheorem{fact}{Fact}[section]
\def\proof{\bigskip \penalty 25\noindent{\bf Proof. }}

\def\endproof{\blackslug \bigskip}

\def\blackslug{\hbox{\hskip 1pt \vrule width 4pt height 6pt depth 1.5pt
  \hskip 1pt}}

\newcommand{\op}[1]{{\mathbf #1}}

\newcommand{\krull}[1]{\left\{ {#1} \right\}}

\newcommand{\bracket}[1]{\left[ {#1} \right]}

\newcommand{\pa}[1]{\left( {#1} \right)}

\newcommand{\abs}[1]{\vert {#1} \vert}

\newcommand{\filtspace}{\pa{\Omega,{\cal F},P, {\mathbf F}}}

\newcommand{\bigfilt}{\{{\cal F}_{t}\}_{t\geq0}}

\newcommand{\Ft}[1]{{\cal F}_{#1}}

\newcommand{\map}[3]{{#1}:{#2}\rightarrow {#3}}

\newcommand{\PA}[2]{{#1}\pa{#2}}

\newcommand{\EpX}[2]{E^{#1}\left[ #2 \right]}

\newcommand{\norm}[2]{\| {#1} \|_{{#2}}}

\newcommand{\half}{\frac{1}{2}}

\newcommand{\dfdx}[2]{\frac{\partial #1}{\partial #2}}

\newcommand{\DfDx}[2]{\frac{d #1}{d #2}}

\newcommand{\ddfdxx}[2]{\frac{{\partial}^{2} #1}{\partial {#2}^2}}

\newtheorem{ex}{Exercise}[section]
 \newcommand{\bex}{\begin{ex}}
 \newcommand{\eex}{\end{ex}}
\begin{document}
\bibliographystyle{acm}
\title{\bf The Pedestrian's Guide to Local Time\thanks{I am very grateful to Mariana Khapko for valuable comments, and for giving me the
necessary motivation to write this paper. Many thanks are also due to Boualem Djehiche for
valuable comments and suggestions.}
}
\author{Tomas Bj{\"{o}}rk \
       \small Department of Finance, \\
       \small Stockholm School of Economics, \\
       \small Box 6501, \\
       \small SE-113 83 Stockholm, \\
       \small SWEDEN \\
       \small tomas.bjork@hhs.se
}
\date{November 23, 2015 \\
\mbox{}\\
{\small Preliminary version\\
\mbox{}\\
Comments are welcome}}

\maketitle

\begin{abstract}
These notes contains an introduction to the theory of Brownian and
diffusion local time, as well as its relations to the Tanaka Formula, the
extended Ito-Tanaka formula for convex functions, the running maximum
process, and the theory of regulated stochastic differential equations. The
main part of the exposition is very pedestrian in the sense that there is a
considerable number of intuitive arguments, including the use of the Dirac
delta function, rather than formal proofs. For completeness sake we have,
however, also added a section where we present the formal theory and give
full proofs of the most important results. In the appendices we briefly
review the necessary stochastic analysis for continuous semimartingales.

\end{abstract}
%\newpage
%\tableofcontents
%\newpage
%%%%%%%%%%%%%%
%  SECTION   %
%%%%%%%%%%%%%%
%%%%%%%%%%%%%%
%  SECTION   %
%%%%%%%%%%%%%%
\newpage
%%%%%%%%%%%%%%%%%%%%%%%%%%%%
\section{Introduction}\label{int} 
%%%%%%%%%%%%%%%%%%%%%
It seems to be a rather well established empirical fact that many students find it
hard to get into the theory of local time. The formal definitions often seem nonintuitive
and the proofs seem to be very technical. As a result, the prospective
student of local time ends up by feeling rather intimidated and simply leaves
the subject.
The purpose of these notes is to show that, perhaps contrary to common
belief, local time is very intuitive, and that half an hour of informal reflection
on the subject will be enough to make the main results quite believable.
The approach of the notes is thus that we take an intuitive view of local
time, we perform some obvious calculations without bothering much about rigor,
we feel completely free to use infinitesimal arguments, and in particular we use
the Dirac delta function.
Arguing in this fashion, we are led to informal proofs for many of the main
results concerning local time, including the relation of Brownian local time to
the running maximum, the Tanaka formula, the extended Ito-Tanaka formula
for convex functions, and the theory of regulated SDE:s.
The hope is that, after having read the intuitive parts of the notes in
Sections 2 - 4, the reader should feel reasonably familiar with the concepts and
central results of local time, and that this should be enough to make him/her
motivated to study the full formal theory.
In order to retain some minimal street cred in mathematical circles I have, 
in Section 5, added full formal proofs of some of the main results. In order
to make the notes more self contained I have also, in the appendices, added a
very brief overview of stochastic analysis for continuous semimartingales. This
includes stochastic dominated convergence and the BDG inequality, as well as
the Kolomogorov continuity criterion.
There is of course no claim to originality in the intuitive approach taken in
the notes. This is probably the standard way of thinking for most people who
are familiar with the theory.
%%%%%%%%%%%%%%%%%%%%%%%%
\section{Local time for deterministic functions} \label{ld}
%%%%%%%%%%%%%%%%%%%%%%%%
In this section we introduce local time for deterministic functions, and derive
some of the most important properties. The reason for spending time on the
rather trivial case of deterministic functions is that a number of concepts and
results, which we later will study in connection with Brownian motion, are
very easily understood within the simpler framework. Lebesgue measure will be
denoted by either $dm(x)$ or $dx$.
%%%%%%%%%%%%%%%%%%%%%%%%
\subsection{Introduction} \label{int}
%%%%%%%%%%%%%%%%%%%%%%%%
Let $X$ be an arbitrary continuous function $\map{X}{R_+}{R}$, with the informal
interpretation that $t \longmapsto X_t$ is a trajectory of a stochastic process. We will sometimes use the notation $X(t)$ instead of $X_t$. We now
consider a Borel set $A \subseteq R$ and we want to study the time spent in the set $A$
by the $X$-trajectory.
\bdf
We define occupation time $T_t(A)$ as the time spent by $X$ in
the set $A$ on the time interval $[0, t]$. Formally this reads as
$$
T_t(A) =\int_0^tI_A \krull{X_s} ds,
$$
where $I_A$ is the indicator function of $A$.
\edf
For a fixed $t$, it is easy to see that $T_t(\cdot)$ is a Borel measure on $R$, and there are
now (at least) two natural questions to ask about this measure.
\ben
\item
Under what conditions is $T_t$ absolutely continuous w.r.t. Lebesgue measure,
i.e. when do we have $T_t(dx) << dx$?
\item
If $T_t(dx) << dm(x)$, we can define $L_t^x$
by
$$
L_t^x= \frac{T_t(dx)}{m(dx)},
$$
and a natural question to ask is under what conditions $L$ is jointly continuous
in $(t, x)$.
\en
The function $L_t$ above is in fact the main actor in the present text, so we write
this as a separate definition.

\bdf
 Assume that $T_t(dx) << dx$. We then define the local time of
$X$ as the Radon-Nikodym derivative
$$
L_t^x= \frac{T_t(dx)}{m(dx)}.
$$
\edf

\begin{remark}\label{vrem}
We note that the physical dimension of local time is time over
distance, i.e. $(velocity)^{-1}$.
\end{remark}
The main object of this text is to study local time for Brownian motion and
diffusion processes, to show that it exists and that it is indeed continuous in
$(t,x)$, to use Brownian local time to obtain an extension of the standard Ito
formula, and to discuss the relation between Brownian local time and the theory
of regulated SDE:s.

%%%%%%%%%%%%%%%%%%%%%%%%%%%%%%%%%%%
\subsection{Absolute continuity of $T_t$ and continuity of $L_t$} \label{ac}
%%%%%%%%%%%%%%%%%%%%%%%%%%%%%%%%%%%%%
In order to get some intuitive understanding of absolute continuity of $T_t(dx)$
and continuity of $L_t^x$
we now study some simple concrete examples.
%%%%%%%%%%%%%%%%%%%%%%%%%%%%%%%%%%%%%%%%%%%
\subsubsection{The case when X is constant} \label{c}
%%%%%%%%%%%%%%%%%%%%%%%%%%%%%%%%%%%%%%%%%%%%%
It is very easy to see that we cannot in general expect to have $T_t(dx) << dm(x)$.
A trivial counterexample is given by any function which is constant on a time
set of positive Lebesgue measure. Consider for example the function $X$ defined
by $X_t=a$ for all $t \geq 0$. For this choice of $X$ we have
$$
T_t(A)=t\cdot I_A\krull{a},
$$
so the measure $T_t(dx)$ corresponds to a point mass of size $t$ at $x=a$. Thus $T_t$ is obviously not
absolutely continuous w.r.t Lebesgue measure.
%%%%%%%%%%%%%%%%%%%%%%%%%%%%%%%%%%%%%%%%%%%
\subsubsection{The case when X is differentiable} \label{d}
%%%%%%%%%%%%%%%%%%%%%%%%%%%%%%%%%%%%%%%%%%%%%
From the previous example we see that in order to guarantee $T_t(dx) << dm(x)$
we cannot allow $X$ to be constant on a time set of positive Lebesgue measure.
It is therefore natural to consider the case when $X$ strictly increasing and for
simplicity we assume that $X$ continuously differentiable. From the definition we
then have
$$
T_t(A)=\int_0^tI_A\krull{X_s}ds.
$$
We now change variable by introducing $x=X_s$, which gives us
\beqarno
s&=&X^{-1}(x),\\
&&\\
ds&=&\frac{dx}{X'\pa{X^{-1}(x)}}.
\eeqarno
Using this we obtain
$$
T_t(A)=\int_A \bracket{\frac{1}{X'\pa{X^{-1}(x)}}I\krull{X_0 \leq x \leq X_t}}dx
$$
From this we see that we do indeed have $T_t(dx) << dm(x)$ and that the local time
is given by
$$
L_t^x=\frac{1}{X'\pa{X^{-1}(x)}}I\krull{X_0 \leq x \leq X_t}
$$
Because of the appearance of the indicator in the expression above we immediately
see that $L$ is not continuous in $t$ and $x$. The fact that we have a term
of the form $1/X'$  is more or less what could be expected from the fact that,
according to Remark \ref{vrem}, local time has dimension $(velocity)^{-1}$. We also see
that we would have a major problem if $X_t'=0$ for some $t$, and this is very
much in line with our discovery in the previous section that we cannot allow
a function which is locally constant if we want to ensure the existence of local
time.

We now go on to consider a more general (not necessarily increasing) function
$X$. We assume that $X$ is continuously differentiable and that the level
sets of $X$ are discrete. In other words we assume that, for every $(t,x)$, the
set $\krull{s\leq t : X_s = x}$ is a finite (or empty) set. This means that we do not
allow functions which are locally constant and we also rule out functions with
oscillatory behavior.
Using more or less the same arguments as for the previous example it is
easy to see that we now have the following result.
\bprop
Assume that $X$ is differentiable with discrete level sets. Then
local time $L_t^x$ exists and is given by the formula
$$
L_t^x=\sum_{X_s=x \atop 0 \leq s \leq t}
\frac{1}{\abs{X_s'}}I\krull{\inf_{s \leq t}X_s \leq x \leq \sup_{s \leq t}X_s}
$$
\eprop

We can immediately note some properties of $L$ for the case of a differentiable $X$ trajectory.
\ben
\item
For fixed $x$, the mapping $t \longmapsto L_t^x$
is non-decreasing.
\item 
As long as $X_t \neq  x$ the mapping $t \longmapsto L_t^x$
is piecewise constant.
\item 
At a hitting time, i.e. when $X_t=x$, then $L^x$  has a jump at $t$ with jump
size 
$\Delta L_t^x=1/\abs{X_t'}$  This is also more or less what could be expected from
Remark \ref{vrem}.
\item 
If $x < inf_s X_s$ or $x > sup_s X_s$ then $L_t^x=0$ for all $t$
\item 
For fixed $t$, the mapping $x \longmapsto L_t^x$ is
continuous except when $x$ is a stationary
point of the $X$ trajectory in which case $L_t^x= + \infty$.
\en

From the discussion above we see that in order to have any chance of obtaining
a local time which is continuous in $(t,x)$ we need to consider a function $X$
which, loosely speaking, has $\abs{X_t'} =+ \infty$  at all points, and which is also rapidly
oscillating. Intuitively speaking we could then hope that the finite sum above
is replaced by an infinite sum with infinitesimal terms. Using more precise
language we thus conjecture that the only $X$ trajectories for which local time
may exist as a continuous function are those with locally infinite variation.
This observation leads us to study local time for Brownian motion and we
will  see that when $X$ is Wiener process, then continuous local time does indeed exist
almost surely.
%%%%%%%%%%%%%%%%%%%%%%%%%%%%%%%%%%%%
\subsection{An integral formula} \label{int}
%%%%%%%%%%%%%%%%%%%%%%%%%%%%%%%%%%%%%%
Let us consider a function $X$ for which local time does exist, but where we make
no assumption of continuity in $t$ or $x$ for $L_t^x$. Recalling the definitions
\beqarno
T_t(A)&=&\int_0^tI_A\krull{X_s}ds,\\
L_t^x&=&\frac{T_t(dx)}{dx}
\eeqarno
we obtain the formula
$$
\int_0^tI_A\krull{X_s}ds=T_t(A)=\int_AL_t^xdx=\int_RI_A(x)L_t^xdx
$$
Using a standard approximation argument we have the following result which
shows that we may replace a time integral with a space integral.
\bprop \label{intprop1}
Let $\map{f}{R}{R}$  be a Borel measurable function in $L^1(L_t^xdx)$
We then have 
$$
\int_0^tf\pa{X_s}ds=\int_Rf(x)L_t^xdx.
$$
\eprop
%%%%%%%%%%%%%%%%%%%%%%%%%%%%%%%%%%%%
\subsection{Connections to the Dirac delta function} \label{dir}
%%%%%%%%%%%%%%%%%%%%%%%%%%%%%%%%%%%%%%
Local time has interesting connections with the Dirac delta function, and we now go on to discuss these.
%%%%%%%%%%%%%%%%%%%%%%%%%%%%%%%%%%%%
\subsubsection{The Dirac delta function} \label{tdir}
%%%%%%%%%%%%%%%%%%%%%%%%%%%%%%%%%%%%%%
In this section we give a heuristic definition of the Dirac function  and discuss
some of its properties.
For a fixed $y \in R$ , the {\bf Dirac delta function} $\delta_y(x)$ can informally be viewed
as a ``generalized function'' defined by the limit
$$
\delta_y(x)=\lim_{\epsilon \rightarrow 0}\frac{1}{2 \epsilon}I\krull{y-\epsilon \leq x \leq y+ \epsilon}
$$
This limit is not a function in the ordinary sense, but we can informally view the expression
$\delta_y(x)dx$ as a unit point mass at $x = y$, so $\delta_y$ is informally the density of the point mass at $y$ w.r.t. Lebesgue measure $dx$. 
We can also view  as a distribution, in
the sense of Schwartz, namely that for any test function $\map{f}{R}{R}$ we
would have
$$
\int_Rf(x)\delta_y(x)dx=f(y).
$$
Using the limit definition above, we also obtain the formulas
\beqarno
\delta_y(x)&=&\delta_0(x-y),\\
\int_Rf(x)\delta_x(y)dx&=&f(y).
\eeqarno
The Dirac function is intimately connected with the Heaviside function $H$
defined by
$$
H_y(x)=
\left\{
\begin{array}{ccc}
0&\mbox{for}& x < y,\\
1&\mbox{for}& x \geq y,
\end{array}
\right.
$$
The Heaviside function is of course not differentiable at $x = y$, but to get some
feeling for what it's ``derivative'' should be, we approximate $H_y$ by
$$
H_y^{\epsilon}(x)=
\left\{
\begin{array}{ccc}
0&\mbox{for}& x < y- \epsilon ,\\
\frac{1}{2 \epsilon} \pa{x-\bracket{y- \epsilon}} &\mbox{for}& y- \epsilon \leq x \leq y+\epsilon , \\
1&\mbox{for}& x \geq y+ \epsilon.
\end{array}
\right.
$$
The derivative of this function is, apart from the points $y \pm \epsilon$, easily seen (draw a figure) to be
given by
$$
\DfDx{H_y^{\epsilon}}{x}(x)=\frac{1}{2 \epsilon} I\krull{y-\epsilon \leq x \leq y+\epsilon}.
$$
As $\epsilon \rightarrow 0$ this is exactly the definition of the Dirac function, and we have given
a heuristic argument for the following result.

\bprop
 The derivative of the Heaviside function is, in the distributional
sense, given by
$$
H_y'(x)=\delta_y(x).
$$
\eprop
%%%%%%%%%%%%%%%%%%%%%%%%%%%%%%%%%%%%
\subsubsection{The delta function and local time} \label{dlt}
%%%%%%%%%%%%%%%%%%%%%%%%%%%%%%%%%%%%%% 
We recall the definition of local time as
$$
L_t^x=\frac{T_t(dx)}{dx}
$$
where occupation time $T_t$ was defined by
$$
T_t(A)=\int_0^tI_A\krull{X_s}ds.
$$
We would thus expect to have
$$
L_t^x=\lim_{\epsilon \downarrow 0}\frac{1}{2 \epsilon} \int_0^t I_{\bracket{x- \epsilon, x + \epsilon}}\krull{X_s}ds
$$
Taking the limit inside the integral we thus expect to have the following result.
\bprop
We may, at least informally, view local time as given by the
expression
$$
L_t^x=\int_0^t\delta_x\pa{X_s}ds,
$$
or on differential form 
$$
dL_t^x=\delta_x\pa{X_t}dt
$$
where the differential  $dL_t^x$ operates on the $t$-variable.
\eprop
%%%%%%%%%%%%%%%%%%%%%%%%%%%%%%%%%%%%
\section{Local time for Brownian motion} \label{lb}
%%%%%%%%%%%%%%%%%%%%%%%%%%%%%%%%%%%%%% 
We now go on to study local time for a Wiener process W.
%%%%%%%%%%%%%%%%%%%%%%%%%%%%%%%%%%%%
\subsection{Basic properties of Brownian local time} \label{lbb}
%%%%%%%%%%%%%%%%%%%%%%%%%%%%%%%%%%%%%% 
Based on the irregular behavior of the Wiener trajectory and the considerations
of Section \ref{ld} we expect occupation time to be absolutely continuous w.r.t.
Lebesgue measure. Based on the informal arguments of Section \ref{d} we also
have some hope that, $L_t^x$
is continuous in $(t,x)$. That this is indeed the case is
guaranteed by a very deep result due to Trotter. We state the Trotter result as
a part of the following theorem where we also collect some results which were
proved already in Section \ref{lb}.
%%%%%%%%%%%%%%%%%%%%%%%%%%%%%%%%%%%%%%%%%%%%%%%
\bth \label{lbbt}
For Brownian motion W, local time $L_t^x$
does exist and the following
hold.
\ben
\item
For almost every $\omega$, the function $L_t^x(\omega )$
is jointly continuous in $(t,x)$.
\item 
$L_t^x$
can be expressed as the limit
$$
L_t^x=\lim_{\epsilon \downarrow 0} \frac{1}{2 \epsilon}\int_0^tI\krull{x- \epsilon \leq W_s \leq x+ \epsilon}ds
$$
\item 
For every fixed $x$ and for almost every $\omega$, the process $t \longmapsto L_t^x(\omega)$ is nonnegative, non decreasing, and increases only on the set $\krull{t \geq 0 : \ W_t(\omega) =x}$.
\item 
Informally we can express local time as
$$
L_t^x=\int_0^t\delta_x\pa{W_s}ds
$$
or on differential form as
$$
dL_t^x=\delta_x\pa{W_t}dt.
$$
\item 
For every bounded Borel function $\map{f}{R}{R}$ we have 
$$
\int_0^tf\pa{W_s}fds= \int_Rf(x)L_t^xdx
$$
\en
\eth
%%%%%%%%%%%%%%%%%%%%%%%%%%%%%%%%%%%%%%%%%%%%%%%%%%%%%%%%%%
From item 3 of this result it is clear that the process $t \longmapsto L_t^x(\omega)$
is a rather
strange animal. Firstly we note that for any fixed $x$, the corresponding level set
$D_x=\krull{t \geq 0:\ W_t=x}$ has Lebesgue measure zero. This follows  from
the fact that occupation time $T_t(x)$ is absolutely continuous w.r.t. Lebesgue
measure, but we can also prove it   by noting that we have $m(D_x)=\int_0^{\infty}I\krull{W_t=x}dt$, 
so we have $\EpX{}{D_x}=\int_0^{\infty}\PA{P}{W_t=x}dt=0$.
Since $m(D_x) \geq 0$ we thus have $m(D_x)=0$ almost surely. One can in fact prove the following
much stronger result.
\bprop
For a fixed $x$, the level set $D_x=\krull{t \geq 0:\ W_t=x}$ is a closed
set of Lebesgue measure zero with the property that every point in $D_x$ is an
accumulation point.
\eprop
The process $L_t^x$ is thus a non decreasing process with continuous trajectories
which increases on a set of Lebesgue measure zero. It is very hard to
imagine what such a trajectory looks like, but one may for example think of the
Cantor function.

We end this section by computing the expected value of local time.
\bprop
Denoting the density of $W_t$ by $p(t, x)$ we have
$$
\EpX{}{L_t^x}=\int_0^tp(s,x)ds.
$$
\eprop
\proof 
Using item 2 of  Theorem \ref{lbbt} and denoting the cdf of $W_t$ by
$F(t, x)$, we have
\beqarno
\EpX{}{L_t^x}&=&\EpX{}{\lim_{\epsilon \downarrow 0} \frac{1}{2 \epsilon}\int_0^tI\krull{x- \epsilon \leq W_s \leq x+ \epsilon}ds}\\
&=&\lim_{\epsilon \downarrow 0} \frac{1}{2 \epsilon}\int_0^t\PA{P}{x- \epsilon \leq W_s \leq x+ \epsilon}ds\\
&=&\lim_{\epsilon \downarrow 0} \frac{1}{2 \epsilon}\int_0^t\krull{F(s, x+ \epsilon) -F(s,x -\epsilon}ds\\
&=&\int_0^tp(s,x)ds.\quad  \endproof
\eeqarno
%%%%%%%%%%%%%%%%%%%%%%%%%%%%%%%%%%%%
\subsection{Local time and Tanaka's formula} \label{tf}
%%%%%%%%%%%%%%%%%%%%%%%%%%%%%%%%%%%%%% 
Let $W$ again be a Wiener process and consider the process $f(W_t)$ where 
$f(x) =\abs{x}$. We would now like to compute the stochastic differential of this process,
and a formal application of the standard Ito formula gives us
$$
df(W_t) = f'(W_t)dW_t + \half f''(W_t)dt.
$$
This formula is of course not quite correct, because of the singularities of the first
and second derivatives of $f$ at $x = 0$. Nevertheless it is tempting to interpret
the derivatives in the distributional sense. We would then have
$$
f'(x)=
\left\{
\begin{array}{rcl}
-1& \mbox{for} &x \leq 0,\\
1 &\mbox{for} &x > 0,
\end{array}
\right.
$$
so $f'$ is a slightly modified Heaviside function. Arguing as in Section \ref{tdir} we
would thus have
$$
f''(x) = 2\delta_0(x)
$$
We now formally insert these expressions for $f'$ and  $f''$ into the Ito formula
above and integrate. We then obtain the expression
$$
\abs{W_t} =\int_0^t \mbox{sgn}(W_s)dW_s+ \int_0^t \delta_0(W_s)ds,
$$
where sgn denotes the sign function defined by
$$
\mbox{sgn}(x)=
\left\{
\begin{array}{rcl}
-1& \mbox{for} &x \leq 0,\\
1 &\mbox{for} &x > 0,
\end{array}
\right.
$$
From Theorem \ref{lbbt} we know that 
$L_t^x=\int_0^t\delta_x(W_s)ds$
so we would formally have
$$
\abs{W_t} =\int_0^t \mbox{sgn} (W_s)dW_s+ L_t^0. 
$$
It is easy to extend this argument to the case when $f(x) = \abs{x - a}$ for some real
number $a$. Arguing as above we have then given a heuristic argument for the
following result.
%%%%%%%%%%%%%%%%%%%%%%%%%%%%%%%%%%%%%%%%%%%%%%%%%
\bth[Tanaka's formula]
For any real number $a$ we have
$$
\abs{W_t-a}=\abs{a}+ \int_0^t \mbox{\em sgn}(W_s-a)dW_s + L_t^a. 
$$
\eth
%%%%%%%%%%%%%%%%%%%%%%%%%%%%%%%%%%%%%%%%%%%%%%
The full formal proof of the Tanaka result is not very difficult. The idea
is simply to approximate $f(x) = \abs{x}$ with a sequence $\krull{f_n}$ of smooth functions
such that $f_n \rightarrow f$, apply the standard Ito formula to each $f_n$, and go to the
limit.

We end this section by connecting the Tanaka formula with the Doob-
Meyer decomposition. We then observe that since $W_t-a$ is a martingale and the
mapping $x \longmapsto \abs{x}$ is convex, the process $\abs{W_t - a}$ is a submartingale. Acccording
to the Doob-Meyer Theorem we can thus decompose $\abs{W_t - a}$ uniquely as
$$
\abs{W_t - a} =\abs{a}+ M_t^a+ A_t^a
$$
where for each $a$, the process $M^a$ is a martingale and $A^a$ is a nondecreasing
process with $M_0^a=A_0^a=0$. Comparing to the Tanaka formula we can thus make the identification
$M_t^a=\int_0^t \mbox{sgn}(W_s-a)dW_s$ and $A_t^a=L_t^a$.
%%%%%%%%%%%%%%%%%%%%%%%%%%%%%%%%%%%%%%%%%%%%%%%%%%%%%%%%%%%%%%%%%%%%%%%%%
\subsection{An alternative definition of local time}
%%%%%%%%%%%%%%%%%%%%%%%%%%%%%%%%%%%%%%%%%%%%%%%%%%%%%%%%%%%%%%%%%%%%%%%
We may in fact use the Tanaka formula to obtain an alternative definition of
Brownian local time.
\bdf
 For any real number $a$, we define local time $L_t^a$ by
$$
L_t^a=\abs{W_t-a}-\abs{a}- \int_0^t\mbox{\em sgn}(W_s-a)dW_s.
$$
\edf
The advantage of this definition is of course that the existence of local
time presents no problem. Continuity in $t$ is obvious, but the hard work is to
show continuity in $(t,a)$, and also to show that we have the property (which
previously was a definition)
$$
T_t(E)=\int_EL_t^ada
$$
for all Borel sets $E$.
%%%%%%%%%%%%%%%%%%%%%%%%%%%%%%%%%%%%%%%%%%%%%%%%%%%
%%%%%%%%%%%%%%%%%%%%%%%%%%%%%%%%%%%%%%%%%%%%%%%%%%%%%%%%%%%%%%%%%%%%%%%%%
\subsection{Regulated Brownian motion} \label{rb}
%%%%%%%%%%%%%%%%%%%%%%%%%%%%%%%%%%%%%%%%%%%%%%%%%%%%%%%%%%%%%%%%%%%%%%%
We now turn to the concept of regulated Brownian motion. This is a topic
which is interesting in its own right, and it is also closely connected to Brownian
local time and to the distribution of the Brownian running maximum. We follow the exposition in \cite{mckean}.
%%%%%%%%%%%%%%%%%%%%%%%%%%%%%%%%%%%%%%%%%%%%%%%%%%%%%%%%%%%%%%%%%%%%%%%%%
\subsubsection{The Skorohod construction} \label{sk}
%%%%%%%%%%%%%%%%%%%%%%%%%%%%%%%%%%%%%%%%%%%%%%%%%%%%%%%%%%%%%%%%%%%%%%%
We start by the following result by Skorohod, which is stated for an arbitrary
continuous function $X$. Loosely speaking we would now like to ``regulate'' $X$ in such a way that
the regulated version of $X$, henceforth denoted by $Z$,  stays positive, and we would also like this regulation to be done in some
minimal way. The following result by Skorohod gives a precise formulation and
solution to the intuitive problem.
%%%%%%%%%%%%%%%%%%%%%%%%%%%%%%%%%%%%%%%%%%%%%%%%%%%%%
\bprop[Skorohod] \label{skt}
%%%%%%%%%%%%%%%%%%%%%%%%%%%%%%%%%%%%%%%%%%%%%%%%
Let $\map{X}{R}{R}$ be an arbitrary continuous function with
$X_0=0$. Then there exists a unique pair of functions $(Z, F)$ such that
\ben
\item
$Z_t = X_t + F_t$, for all $t \geq 0$.
\item 
$Z_t \geq 0$, for all $t \geq 0$.
\item 
$F$ is nondecreasing with $F_0=0$.
\item 
$F$ increases only when $Z_t=0$.
\en
\eprop
\proof We see immediately that if a pair $(Z, F)$ exists, then $F_t \geq X_t^-$ for all $t$. 
and since $F$ is nondecreasing we must in fact have $F_t \geq \sup_{s \leq t}X_s^-$.
This leads us to conjecture that in fact $F_t = \sup_{s \leq t}X_s^-$ and it is easy to see that with
this choice of $F$ all the requirements of the proposition are satisfied. To prove
uniqueness, let us assume that there are two solutions $(Z, F)$ and $(Y,G)$ to the
Skorohod problem so that
\beqarno
Z_t&=&X_t + F_t,\\
Y_t&=&X_t + G_t.
\eeqarno
Define $\xi$ by $\xi_t = Z_t-Y_t=F_t-G_T$.. Using the fact that $\xi$ is continuous and of
bounded variation then have
$$
d\xi_t^2=2\xi_td\xi_t=2(Z_t-Y_t)(dF_t-dG_t).
$$
Since $F$ is constant off the set $\krull{Z_t = 0}$ we have $Z_tdF_T=0$ and in the same way
we have $Y_tdG_t=0$. We thus obtain
$$
\xi_t^2=-2\int_0^t\pa{Z_sdG_s + Y_sdF_s}\leq 0,
$$
but since $\xi^2$ obviously is positive we conclude that $\xi =0$.
\endproof 

We note that condition that $F$ increases only on the set $\krull{Z_t = 0}$ formalizes
the idea that we are exerting a minimal amount of regulation.
%%%%%%%%%%%%%%%%%%%%%%%%%%%%%%%%%%%%%%%%%%%%%%%%%%%%%%%%%%%%%%%%%%%%%%%%%
\subsection{Local time and the running maximum of $W$} \label{rm}
%%%%%%%%%%%%%%%%%%%%%%%%%%%%%%%%%%%%%%%%%%%%%%%%%%%%%%%%%%%%%%%%%%%%%%%
In this section we will investigate the relations between local time, regulated
Brownian motion, and the running maximum of a Wiener process. We start by
recalling the Tanaka formula
$$
\abs{W_t}=\int_0^t\mbox{sgn}(W_s)dW_s+L_t^0.
$$
where $L^0=L^0(W)$ is local time at zero for $W$. We note that the process $M$
defined by
$$
M_t=\int_0^t\mbox{sgn}(W_s)dW_s
$$
is a martingale and, since we obviously have $(dM_t)^2 = dt$, it follows from the
Levy characterization (see Corollary A.1) that $M$ is in fact a Wiener process,
so we can write $M = ˜\tilde{W}$. We may thus write the Tanaka formula as
$$
\abs{W_t} =\tilde{W}_t+L_t^0(W).
$$
We obviously have $\abs{W_t}\geq 0$ and we we recall that $L^0$ only increases on the set
$\krull{W_t = 0}$ i.e. on the set $\krull{\abs{W_t} = 0}$. Let us now apply the Skorohod regulator
to the process ˜$\tilde{W}$. We then obtain
$$
Z_t = ˜\tilde{W}_t + F_t,
$$
where $Z \geq 0$, and $F$ is increasing only on $\krull{Z_t = 0}$. Comparing these two
expressions we see that the pair $\pa{\abs{W},L^0}$ is in fact the solution to the Skorohod
regulator problem for the process $\tilde{W}$ so we have $Z_t=\abs{W_t}$ and $F_t=L_t^0(W)$.
Furthermore, we know from the Skorohod theory that $F_t=\sup_{s \leq t}\tilde{W}_s^-$,
so we obtain
\beq \label{rm8}
L_t^0(W)=\sup_{s \leq t}\tilde{W}_s^-,
\eeq
so in particular we have the formula
\beq \label{rm10}
\abs{W_t}=\tilde{W}_t+\sup_{s \leq t}\tilde{W}_s^-
\eeq
which we will need later on. 

Going back to formula  \cref{rm8}
we emphasize that $L_t^0(W)$ is local time for $W$ but that the supremum in the
right hand side of this equality is for $\tilde{W}$, so we are not allowed to identify $L^0(W)$
with $\sup_{s \leq t}{W}_s^-$.
We note however that we have
$$
\tilde{W} \stackrel{d}{=}W
$$
where $\stackrel{d}{=}$ denotes equality in distribution. Using this and the symmetry of the
Wiener process we obtain
$$
L_t^0(W)=\sup_{s \leq t}\tilde{W}_s^- \stackrel{d}{=} \sup_{s \leq t}{W}_s^- \stackrel{d}{=} \sup_{s \leq t}{W}_s.
$$
We have thus proved the following result.
%%%%%%%%%%%%%%%%%%%%%%%%%%%%%%%%%%%%%%%%%%%%
\bprop
%%%%%%%%%%%%%%%%%%%%%%%%%%%%%%%%%
Let $W$ be a Wiener process and define the running maximum
by
$$
S_t = \sup_{s \leq t}{W}_s.
$$
We then have the relation
$$
L_t^0(W)\stackrel{d}{=} S_t.
$$
\eprop
%%%%%%%%%%%%%%%%%%%%%%%%%%%%%%%%%%%%%%%%
We may in fact improve this result considerably. From \cref{rm10} we have the relation
$$
\abs{{W}_t} = \tilde{W}_t+ \sup_{s \leq t} \tilde{W}_s^-
$$
Let us now perform a simple change of notation. We denote $W$ by $\hat{W}$, and we denote $\tilde{W}$ by $-W$. Noting that $W$ is a Wiener process we have  the formula
$$
\abs{\hat{W}_t} = -W_t+ \sup_{s \leq t} {\pa{-W_s}^-}.
$$
We can thus view this formula as  the Kolmogorov regulator applied to the Wiener process $-W$.
We also have
$$
\sup_{s \leq t} {\pa{-W_s}^-}=\sup_{s \leq t}\pa{W_s}=S_t,
$$
so
$$
\abs{\hat{W}_t}=-W_t+S_t.
$$
On the other hand, we have the trivial equality
$$
S_t -W_t = -W_t + S_t.
$$
so we obtain
$$
(S -W, S) = \pa{\abs{\hat{W}},S}
$$
From Skorohod and the Tanaka formula we know that $S$ is local time for $\hat{W}$ so
$S_t =L_t^0(\hat{W})$. This gives us
$$
(S -W, S) =\pa{\abs{\hat{W}},L^0(\hat{W})}
$$
and, since $W$ and $\hat{W}$ have the same distribution, we have thus proved the following
result.
%%%%%%%%%%%%%%%%%%%%%%%%%%%%%%%%%%%%%%%%%%%%
\bprop
%%%%%%%%%%%%%%%%%%%%%%%%%%%%%%%%%%%%%%%%%%
Let $W$ be a Wiener process, $S$ its running maximum, and $L^0$
it's local time at zero. We then have the  distributional equality
$$
(S -W, S) \stackrel{d}{=}\pa{\abs{{W}},L^0({W})}.
$$
\eprop
%%%%%%%%%%%%%%%%%%%%%%%%%%%%%%%%%%%%%%%%%
%%%%%%%%%%%%%%%%%%%%%%%%%%%%%%%%%%%%%%%%%%%%%%%%%%%%%%%%%%%%%%%%%%%%%%%%%
\section{Local time for continuous semimartingales} \label{ls}
%%%%%%%%%%%%%%%%%%%%%%%%%%%%%%%%%%%%%%%%%%%%%%%%%%%%%%%%%%%%%%%%%%%%%%%
In this section we extend our analysis of local time to the case of a continuous
semimartingale $X$ with decomposition
$$
X_t = A_t +M_t,
$$
where $A$ is a continuous process of local finite variation and $M$ is a continuous local martingale (see Appendix \ref{as}). 
The most common special case of this is of course when $X$ is an Ito process
and thus has dynamics of the form
$$
dX_t = \mu_tdt +\sigma_tdW_t
$$
where $W$ is Wiener.
%%%%%%%%%%%%%%%%%%%%%%%%%%%%%%%%%%%%%%%%%%%%%%%%%%%%%%%%%%%%%%%%%%%%%%%%%
\subsection{Definition, existence, and basic properties} \label{lsd}
%%%%%%%%%%%%%%%%%%%%%%%%%%%%%%%%%%%%%%%%%%%%%%%%%%%%%%%%%%%%%%%%%%%%%%%
It would seem natural to define occupation time $T_t(A)$ in the usual way as
$T_t(A)=\int_0^tI_A(X_s)ds$ and then define local time as $L_t^x=T_t(dx)/dm$.
It is perfectly possible to develop a theory based on these definitions, but instead it
has become standard to measure occupation time with respect to the quadratic
variation process $\qvar{X}$. The concept of quadratic variation belongs to general
semimartingale theory, but for Ito processes it is very simple. See Appendix \ref{qv}
for details.
%%%%%%%%%%%%%%%%%%%%%%%%%%%%%%%%%%
\bdf
For an Ito process with dynamics
$$
dX_t = \mu_tdt +\sigma_tdWt
$$
we define the quadratic variation process $\qvar{X}$ by
$$
\qvar{X}_t=\int_0^t\sigma_s^2ds.
$$
We can also write this as
$$
d\qvar{X}_t=\sigma_t^2dt, \quad \mbox{or as}\quad  d\qvar{X}_t=\pa{dX_t}^2
$$
\edf
%%%%%%%%%%%%%%%%%%%%%%%%%%%%%%%%%%%%%%%%%%
We may now define (scaled) occupation time.
%%%%%%%%%%%%%%%%%%%%%%%%%%%%%%%%%%%%%%%%%%
\bdf
For any Borel set $A \subseteq R$ we define occupation time $T_t(A)$ by
$$
T_t(A)=\int_0^tI_A(X_s)d\qvar{X}_s,
$$
so for an Ito process we have
$$
T_t(A)=\int_0^tI_A(X_s)\sigma_s^2ds.
$$
\edf
%%%%%%%%%%%%%%%%%%%%%%%%%%%%%%% 
We note that for a standard Wiener process $W$, we have $d\qvar{W}_t=dt$, so for
a Wiener process the new definition coincides with the old one. Local time is
now defined as before.
%%%%%%%%%%%%%%%%%%%%%%%%%%%%%%%%%%%%
\bdf
If $T_t(dx) << dm$, then we define local time $L_t^x$ by
$$
L_t^x=\frac{T_t(dx)}{m(dx)}
$$
\edf
%%%%%%%%%%%%%%%%%%%%%%%%%%%%%%%%%%%%
With the new definition of occupation time we also see that that if $\sigma \equiv 0$ (so
$X$ is of bounded variation), then $T_t \equiv 0$. Thus local time does (trivially) exist
and we have $L \equiv 0$

The interesting case is of course when $\sigma \neq 0$. We then have the following result,
where the first item is deep and the other items are more or less expected.
%%%%%%%%%%%%%%%%%%%%%%%%%%%%%%%%%%
\bth \label{ltdt}
 Assume that $X$ is a continuous semimartingale with canonical
decomposition $X = A +M$. Then the following hold.
\ben
\item
Local time $L_t^x$ exists. Furthermore, there is a version of $L$ which is continuous
in $t$ and cadlag in $x$.
\item
The jump size of $L$ in the $x$ variable is given by
$$
L_t^x-L_{t}^{x-}=2\int_0^tI\krull{X_s=x}dA_s.
$$
\item
If $X$ is a local martingale, then $L$ is almost surely continuous in $(t, x)$.
\item
For every fixed $x$, the process $ t \longmapsto L_t^x$
is nondecreasing and only increases on the
set $\krull{t \geq 0:\ X_t = x}$.
\item
For every Borel set $A$ we have
$$
\int_0^tI_A(X_s)d\qvar{X}_s=\int_RI_A(x)L_t^xdx.
$$
\item
For every bounded Borel function $\map{f}{R}{R}$ we have
$$
\int_0^tf(X_s)d\qvar{X}_s=\int_Rf(x)L_t^xdx.
$$
\item
We have the informal interpretation
$$
L_t^x=\int_0^t\delta_x(X_s)d\qvar{X}_s.
$$
\en
\eth
%%%%%%%%%%%%%%%%%%%%%%%%%%%%%%%%%%%%%%%%%%%%%%
Note that a discontinuity in $x \longmapsto L_t^x$  
can only occur if the finite variation process
$A$ charges the set $\krull{X_s = x}$. A concrete example is given in Section \ref{anex}.
%%%%%%%%%%%%%%%%%%%%%%%%%%%%%%%%%%%%%%%%%%%%%%
%%%%%%%%%%%%%%%%%%%%%%%%%%%%%%%%%%%%%%%%%%%%%%%%%%%%%%%%%%%%%%%%%%%%%%%%%
\subsection{A random time change} \label{rt}
%%%%%%%%%%%%%%%%%%%%%%%%%%%%%%%%%%%%%%%%%%%%%%%%%%%%%%%%%%%%%%%%%%%%%%%
We now prove that local time is preserved under a random time change (see Appendix
\ref{rt} for definitions and details on random time changes). Let us therefore
consider a semimartingale $X$, and a smooth random change of time $t \longmapsto C_t$. We
define the process $X^C$ by
$$
X_t^C=X_{C_t},
$$
and we denote the local time of $X$ and $X^C$  by $L(X)$ and $L(X^C)$ respectively.
%%%%%%%%%%%%%%%
\bprop \label{rtt}
With notation as above we have
$$
L_t^a(X^C)=L_{C_t}^a(X)
$$
or equivalently
$$
L(X^C)=L(X)^C.
$$
\eprop
%%%%%%%%%%%%%%%%
\proof
By definition we have
$$
L_t^a(X)=\int_0^t \delta_a(X_s)d\qvar{X}_s,
$$
so we have
$$
L_{C_t}^a(X)=\int_0^{C_t} \delta_a(X_s)d\qvar{X}_s,
$$
from which it follows that
\beq \label{rtt5}
dL_{C_t}^a(X)=\delta_a (X_{C_t})d\qvar{X}_{C_t}.
\eeq
On the other hand, we have by definition
$$
L_t^a(X^C)=\int_0^t\delta_a(X_s^C)d\qvar{X^C}_s,
$$
so
$$
L_t^a(X^C)=\delta_a(X_{C_t})d\qvar{X^C}_t.
$$
From Proposition \ref{rt5} we have $\qvar{X^C}_t=\qvar{X}_{C_t}$ 
so we have
$$
dL_t^a(X^C)=\delta_a(X_{C_t})d\qvar{X}_{C_t}.
$$
Comparing this to \cref{rtt5} we thus have
$$
dL_t^a(X^C)=dL_{C_t}^a(X).\quad
\endproof
$$
%%%%%%%%%%%%%%%%%%%%%%%%%%%%%%%%%%%%%%%%%%%%%%%%%%%%%%%%%%%%%%%%%%%%%%%%%
\subsection{The Tanaka formula} \label{rta}
%%%%%%%%%%%%%%%%%%%%%%%%%%%%%%%%%%%%%%%%%%%%%%%%%%%%%%%%%%%%%%%%%%%%%%%
Assume again that $X$ is a semimartingale and consider the process $Z$ defined
by
$$
Z_t=\abs{X_t}.
$$
We now apply the Ito formula exactly like in the case of Brownian motion.
Interpreting the derivatives of $\abs{x}$ in the distributional sense we would then
formally obtain
$$
\abs{X_t}=\abs{X_0}+ \int_0^t\mbox{sgn}(X_s)dX_s + \int_0^t \delta_0(X_s)d\qvar{X}_s
$$
and from item 7 of Theorem \ref{ltdt}, we conclude that we can write this as
$$
\abs{X_t}=\abs{X_0}+ \int_0^t\mbox{sgn}(X_s)dX_s + L_t^0(X).
$$
Applying the same argument to the function $\abs{x - a}$ we obtain, at least informally,
the following result.
%%%%%%%%%%%%%%%%%%%%%%%%%%%%%%%%%%%%%%%
\bth[Tanaka's formula]
%%%%%%%%%%%%%%%%%%%%%%%%%%%%%%%%%%%%%%%
 With notation as above we have, for every
fixed $a \in R$,
$$
\abs{X_t -a}=\abs{X_0-a}+ \int_0^t\mbox{\em sgn}(X_s-a)dX_s + L_t^a(X).
$$
\eth
%%%%%%%%%%%%%%%%%%%%%%%%%%%%%%%%%%%%%%%%%%%%%%%%
The Tanaka formula leads us to an alternative definition of local time. This
alternative way of defining local time for a diffusion is in fact the most common
definition in the literature.
%%%%%%%%%%%%%%%%%%%%%%%%%%%%%%%%%%%
\bdf
Assume that $X$ is a continuous semimartingale. We may then define local time $L_t^a(X)$ by the formula
$$
L_t^a(X)=\abs{X_t -a}-\abs{X_0-a}- \int_0^t\mbox{\em sgn}(X_s-a)dX_s.
$$
\edf
%%%%%%%%%%%%%%%%%%%%%%%%%%%
%%%%%%%%%%%%%%%%%%%%%%%%%%%%%%%%%%%%%%%%%%%%%%%%%%%%%%%%%%%%%%%%%%%%%%%%%
\subsection{An extended Ito formula} \label{exi}
%%%%%%%%%%%%%%%%%%%%%%%%%%%%%%%%%%%%%%%%%%%%%%%%%%%%%%%%%%%%%%%%%%%%%%%
Let us consider a convex function $\map{f}{R}{R}$, so $f'$ exists almost everywhere. It
furthermore follows from the convexity of $f$ that we may interpret the (distributional) second
derivative $f''$ as a positive measure $f''(dx)$. In most practical applications this
means that we can view $f''$ as consisting of two parts: an absolutely continuous
part, and a finite number of point masses. In most cases we can therefore write
$f''(dx)$ as
$$
f''(dx)=g(x)dx+ \sum_ic_i \delta_{a_i}(x)dx
$$
where the density $g$ is the absolutely continuous part of $f''$, and the point mass
at $x = a_i$ has size $c_i$. In the special case that $f \in C^2$ we would  thus have 
$$
f''(dx)=\DDfDxx{f}{x}(x)dx
$$
If $X$ is a continuous semimartingale, a formal application of the Ito formula gives us
$$
f(X_t)=f(X_0)+ \int_0^tf'(X_s)dX_s + \half \int_0^t f''(X_s)d\qvar{X}_s,
$$
and from Theorem \ref{ltdt}, we conclude that we can write the last term as
$$
\int_0^t f''(X_s)d\qvar{X}_s=\int_0^tL_t^af''(a)da.
$$
The fact that we can view $f''$ as a measure finally allows us to write
$$
\int_0^t f''(X_s)d\qvar{X}_s=\int_0^tL_t^af''(da).
$$
We have thus given a heuristic argument for the following extended Ito formula
which, among other things, streamlines the argument behind the Tanaka
formula.
%%%%%%%%%%%%%%%%%%%%%%%%%%%%%
\bth[Extended Ito formula]
Assume that $X$ is a continuous semimartingale
and that $f$ is the difference between two convex functions. Then we
have
$$
f(X_t)=f(X_0)+ \int_0^tf_{-}'(X_s)dX_s + \half \int_0^t f''(X_s)d\qvar{X}_s,
$$
where $f_{-}'$ is the left hand derivative (which exists everywhere) of $f$, and $f''$ is
interpreted in the distributional sense. Alternatively we may write
$$
f(X_t)=f(X_0)+ \int_0^tf_{-}'(X_s)dX_s + \half \int_0^tL_t^af''(da).
$$
\eth
%%%%%%%%%%%%%%%%%%%%%%%%%%%%%%%%%%%
%%%%%%%%%%%%%%%%%%%%%%%%%%%%%%%%%%%%%%%%%%%%%%%%%%%%%%%%%%%%%%%%%%%%%%%%%
\subsection{An example: The local time at zero of $\abs{W}$} \label{anex}
%%%%%%%%%%%%%%%%%%%%%%%%%%%%%%%%%%%%%%%%%%%%%%%%%%%%%%%%%%%%%%%%%%%%%%%
To see that a continuous semimartingale may indeed have a local time which is non-continuous in the $x$-variable, it is
instructive to study the local time of the process
$$
X_t=\abs{W_t}.
$$
Using the reflection principle of Brownian motion, one would
guess that
$$
L_t^x(X)=2L_t^x(W), \quad \mbox{for $x \geq 0$},
$$
whereas
$$
L_t^x(X)=0,\quad \mbox{for $x < 0$}.
$$
This conjecture is in fact correct, and to see this formally, let us in particular study $L_t^0(X)$.
From Tanaka we have
$$
dX_t = \mbox{sgn}(W_t)dW_t + dL_t^0(X).
$$
Again from the Tanaka formula, and the fact that $X_t=\abs{X_t}$, we have
$$
dX_t=\mbox{sgn}(X_t)dX_t + dL_t^0(X)=\bracket{1-2I\krull{X_t=0}}dX_t+ dL_t^0(X).
$$
We thus obtain
\beqarno
dX_t&=&dX_t-2I\krull{X_t=0}\mbox{sgn}(W_t)dW_t-2I\krull{X_t=0}dL_t^0(W) + dL_t^0(X)\\
&=&dX_t-2L_t^0(W)+dL_t^0(X),
\eeqarno
where we have used the facts that $I\krull{X_t = 0}\mbox{sgn}(W_t)dW_t = 0$, and that we
have $I\krull{X_t = 0}dL_t^0(W)=dL_t^0(W)$.  We have thus proved (as expected) that
$$
L_t^0\pa{\abs{W}}=2L_t^0(W).
$$
Since we obviously have $L_t^x(X)=0$ for $x < 0$, we thus see that $L_t^x(X)$ has a
discontinuity at $x=0$ of size
$$
L_t^0(X)-L_t^{0-}(X)=2L_t^0(W).
$$
and this is exactly what is stated in Theorem \ref{ltdt}.

%%%%%%%%%%%%%%%%%%%%%%%%%%%%%%%%%%%%%%%%%%%%%%%%%%%%%%%%%%%%%%%%%%%%%%%%%
\subsection{Regulated SDE:s} \label{rs}
%%%%%%%%%%%%%%%%%%%%%%%%%%%%%%%%%%%%%%%%%%%%%%%%%%%%%%%%%%%%%%%%%%%%%%%
In this section we will extend the Skorohod construction of regulated Brownian
motion in Section \ref{rb} to the case of a regulated diffusion. The basic idea is the
following.
\bei
\item
We consider an SDE of the form
\beqarno
dX_t &=& \mu(X_t)dt +\sigma(X_t)dW_t,\\
X_0 &=& 0
\eeqarno
\item
We would now like to ``regulate'' the process $X$ in such a way that it stays
positive, and we want to achieve this with a minimal (in some sense) effort.
\item
We thus add an increasing process $F$ to the right hand side of the SDE
above, such that $F$ only increases when $X_t = 0$.
\ei
We start by formally defining the concept of a regulated SDE.
%%%%%%%%%%%%%%%%%%%%%%%%%%%%%%%
\bdf
For given functions $\mu$ and $\sigma$, consider the regulated SDE
\beqarno
dX_t &=& \mu(X_t)dt +\sigma(X_t)dW_t +F_t,\\
X_0 &=& 0
\eeqarno
We say that a pair of adapted processes $(X, F)$ is a solution of the regulated
SDE above if the following conditions are satisfied
\bei
\item
The pair $(X, F)$ does indeed satisfy the regulated SDE above.
\item
For all $t$ we have $X_t \geq 0$,  $P-a.s$.
\item
The process $F$ is increasing and $F_0=0$.
\item
$F$ is  increasing only on the set $\krull{X_t=0}$.
\ei
\edf
We now have the following result by Skorohod.
%%%%%%%%%%%%%%%%%%%%%%%%%%%%%
\bth[Skorohod] 
%%%%%%%%%%%%%%%%%%%%%%%%%%%%
Assume that $\mu$ and $\sigma$ are $C^1$ and that they
satisfy the Lipshitz conditions
$$
\abs{\mu (x) - \mu (y)}\leq K\abs{x-y},\quad \abs{\sigma (x) - \sigma (y)}\leq K\abs{x-y}
$$
for some constant $K$. Then the regulated SDE
\beqarno
dX_t &=& \mu(X_t)dt +\sigma(X_t)dW_t +F_t,\\
X_0 &=& 0
\eeqarno
has a unique solution $(X,F)$. We can furthermore identify $F$ with local time at
zero of $X$, so
$$
F_t=L_t^0(X).
$$
\eth

\proof  We divide the proof, where we follow \cite{mckean}, into several steps.\\
\mbox{}\\
%%%%%%%%%%%%%%%%%%%%%%%%%%%%%%%%%%%%%
{\bf 1. Uniqueness.}\\
Consider two solutions $(X,F)$ and $(Y,G)$. The difference $Z=X-Y$ will then
satisfy the equation
$$
dZ_t=Z_t \hat{\mu}_tdt + Z_t\hat{\sigma}_tdW_t + dF_t - dG_t
$$
where
$$
\mu_t=
\left\{
\begin{array}{ccl}
\frac{\mu(X_t)-\mu(Y_t)} {X_t-Y_t}& \mbox{when} &X_t \neq Y_y,\\
&&\\
\mu'(X_t) &\mbox{when} &X_t=Y_t,
\end{array}
\right.
$$
and similarly for $\hat{\sigma}_t$.
From Ito we obtain
\beqarno
dZ_t^2&=&2Z_tdZ_t + \pa{dZ_t}^2\\
&=&Z_t^2\krull{2 \hat{\mu}_t + \hat{\sigma}_t^2}dt + 2Z_t^2\hat{\sigma}_tdW_t + 2(X_t-Y_t)\krull{dF_t-dG_t}
\eeqarno
We easily see that $(X_t-Y_t)\krull{dF_t-dG_t}=-\pa{Y_tdF_t+X_tdG_t}\leq 0$. Furthermore,
the Lipschitz conditions on $\mu$ and $\sigma$ imply that $\hat{\mu}$ and $\hat{\sigma}$ are bounded by
$K$. Taking expectations and defining $h$ by $h_t=\EpX{}{Z_t^2}$ we thus obtain
$$
\DfDx{h}{t} \leq h_t \pa{2K+ K^2}, \quad h_0=0.
$$
so Gronwalls inequality gives us $h\equiv 0$ and $Z \equiv 0$.\\
\mbox{}\\
%%%%%%%%%%%%%%%%%%%%%%%%%%%%%%%%%%%%%
{\bf 2. Existence for the case $\mu =0$, $\sigma =1$}\\
For this case the regulated SDE takes the form
$$
dX_t=dW_t + dF_t 
$$
so this case is already covered by Proposition \ref{skt}.\\
\mbox{}\\
%%%%%%%%%%%%%%%%%%%%%%%%%%%%%%%%%%%%%
{\bf 3. Existence for the case when $\mu =0$ and $\sigma \neq 0$.}\\
We now have an equation of the form
$$
dX_t=\sigma (X_t)dW_t + F_t.
$$
From step 2 we know that there exists a solution $(X^0,F^0)$ to the equation
$$
dX_t^0=dW_t^0 + F_t^0,
$$
where $W^0$ is a Wiener process. We now define the random time change $t \longmapsto C_t$
by
$$
t=\int_0^{C_t}\frac{1}{\sigma^2(X_s)}ds
$$
and define$X$ and $F$ by
$$
X_t=X_{C_t}^0, \quad F_t = F_{C_t}^0.
$$
It now follows from Proposition \ref{rt23} that $(X,F)$ is the solution to
$$
dX_t=\sigma (X_t)dW_t + dF_t 
$$
where $W$ is a Wiener process.\\
\mbox{}\\
%%%%%%%%%%%%%%%%%%%%%%%%%%%%%%%%%%%%%
{\bf 4. Existence for the case when $\mu \neq 0$ and $\sigma \neq 0$.}\\
This case is easily reduced to step 3 by a Girsanov transformation, and we
note that the properties of $F$ are not changed under an absolutely continuous
measure transformation.\\
\mbox{}\\
%%%%%%%%%%%%%%%%%%%%%%%%%%%%%%%%%%%%%
{\bf 4. Identifying $F$ as local time}\\ 
In step 2 above we know from Section \ref{rm} that we have the identification
$F_t^0=L_t^0(X)$. When we perform the time change in step 3 the local time property
is preserved due to Proposition \ref{rtt}. Finally, the Girsanov transformation in
step 4 does not affect this identification (although it will of course change the
distribution of $L$).
\endproof
%%%%%%%%%%%%%%%%%%%%%%%%%%%%%%%%%%%%%%%%%%%%%%%%%%%%%%%%%%%%%%%%%%%%%%%%%
\section{Some formal proofs} \label{fp}
%%%%%%%%%%%%%%%%%%%%%%%%%%%%%%%%%%%%%%%%%%%%%%%%%%%%%%%%%%%%%%%%%%%%%%%
In this section we present (at least partial) formal proofs for the main results
given above, and for most of the proofs we follow \cite{kall}. These proofs require some
more advanced techniques and results from stochastic analysis. We start with
the definition of local time.
%%%%%%%%%%%%%%%%%%%%%%%%%%
\bdf \label{fp1} Let $X$ be a continuous semimartingale with canoncial decomposition 
$X =M + A$. The local time $L_t^x$
is defined by
$$
L_t^x=\abs{X_t}-\abs{X_0} - \int_0^t\mbox{\em sgn}(X_t-x)dX_t.
$$
\edf
%%%%%%%%%%%%%%%%%%%%%%%%%%%%%%%%%%%%%%%%%%%%%%%%%%%%%%%%%%%%%%%%%%%%%%%%%
\subsection{Basic properties of $L$} \label{bp}
%%%%%%%%%%%%%%%%%%%%%%%%%%%%%%%%%%%%%%%%%%%%%%%%%%%%%%%%%%%%%%%%%%%%%%%
The first result is surprisingly easy to prove.
%%%%%%%%%%%%%%%%%%% 
\bprop 
Local time has the following properties
\bei
\item
For every fixed $x$, the process $t \longmapsto L_t^x$ is nonnegative, continuous and nondecreasing.
\item
For every fixed $X$, the process $L^x$ is supported by the set $\krull{Xt = x}$ in the
sense that 
$$
\int_0^{\infty}I\krull{X_t \neq x}dL_t^x=0
$$
where the differential $dL_t^x$ operates in the $t$ variable.
\item
We have the representation
$$
L_t^x=-\inf{s \leq t}\int_0^s \mbox{\em sgn}(X_u-x)dX_u.
$$
\ei
\eprop
%%%%%%%%%%%%%%
\proof 
We restrict ourselves to the case when $X$ is a local martingale. It is
enough to do the proof for the case $x=0$. For each $h >0$ we can find a convex
$C^2$-function $f_h$ such that $f_h(x)=\abs{x}$ for $ x \leq 0$ and $f_h(x)=x-h$ for $x \geq h$. We
note that
\beqarno
f_h''(x)&=&0 \quad \mbox{for $x \leq 0$ and $x \geq h$},\\
f_h(x)&\rightarrow& \abs{x},\quad \mbox{as $h \rightarrow 0$},\\
f_h'(x)&\rightarrow& \mbox{sgn}(x), \quad \mbox{as $h \rightarrow 0$}.
\eeqarno
We now define the process $Z_t^h$ by
$$
Z_t^h=f_h(X_t) - \int_0^tf_h'(X_s)dX_s.
$$
Since $f_h(x)$ is a a good approximation of $\abs{x}$, we expect that 
$Z_t^h$ should be a good approximation to $L_t^0$.
We obviously have $f_h(X_T) \rightarrow \abs{X_t}$. Furthermore, since
$f_h'(X_t) \rightarrow \mbox{sgn}(X_t)$ and $\abs{f_h'(X_t)}\leq 1$ we conclude from the stochastic dominated
convergence Theorem \ref{dc} that
$$
\pa{Z^h-L^0}^{\star} \stackrel{P}{\rightarrow} 0
$$
where for any process $Y$ we define $Y^{\star}$ by $Y_t^{\star}= \sup_{s \leq t} \abs{Y_s}.$
We have thus seen that $Z_t^h$
converges uniformly in probability to $L_t^0$.
We now apply the Ito formula to $Z_t^h$ and obtain
$$
Z_t^h=\int_0^tf_h''(X_s)d\qvar{X}_s.
$$
Since $f_h$ is convex and $\qvar{X}$ is nondecreasing we see that $Z_t^h$
is increasing, and since $Z_t^h$ converges to $L_t^0$ we conclude that $L_t^0$
is nondecreasing. It is also obvious from the definition that $L_t^0$
is continuous. We have thus proved the first statement in the proposition.

To show the second statement we note that since $f_h''=0$ outside $\bracket{0,h}$ we
trivially have
$$
\int_0^tI\krull{X_s \notin \bracket{0,h}}dZ_s^h=\int_0^tI\krull{X_s \notin \bracket{0,h}}f_h''(X_s)d\qvar{X}_s =0.
$$
Going to the limit this gives us
$$
\int_0^tI\krull{X_s \neq 0}dL_s^0=0.
$$

The third statement follows from the Skorohod result. \endproof

We now go on to prove continuity in $(t,x)$ of $L_t^x$.
The obvious idea is to
use the Kolmogorov criterion (see Theorem \ref{kc}) for continuity, and this will
require the Burkholder-Davis-Gundy inequality \cref{qv2}.
%%%%%%%%%%%%%%%%%%%%%%%%%%%%%%%%%%%%
\bprop[Continuity of $L_t^x$]
With $L_t^x$ defined by Definition \ref{fp1} the following hold.
\bei
\item
There is a version of $L$ which is continuous in $t$ and cadlag in $x$.
\item
The jump size of $L_t$ in the $x$ variable is given by
\beq \label{fp2}
L_t^x-L_t^{x-}=2\int_0^tI\krull{X_s=x}dA_s.
\eeq
\item
If $X$ is a local martingale, then $L$ is almost surely continuous in $(t,x)$.
\ei
\eprop
%%%%%%%%%%%%%%%%%%%
\proof
By definition and with $M_t=\int_0^t \sigma_sdW_s$ we have
$$
L_t^x=\abs{X_t-x}-\abs{X_0-x}- \int_0^t\mbox{sgn} (X_s-x)dM_s - \int_0^t \mbox{sgn} (X_s-x)dA_s.
$$
The term $\abs{X_t-x}$ is obviously continuous in $(t,x)$ and the jump size of the
$dA$-integral is clearly given by the right hand side of \cref{fp2}. It thus remains to
prove continuity of the integral term
$$
I_t^x=\int_0^t\mbox{sgn} (X_s-x)dM_s.
$$
By localization we may assume that the processes $X_t-X_0$, ${\qvar{X}_t^{\half}}$, and $\int_0^t \abs{dA_s}$
are bounded by some constant $C$. For any $ x < y$ the BDG inequality \cref{qv2} gives us
\beqarno
\EpX{}{\pa{I^x-I^y}_t^{\star p}}&=&2^p\EpX{}{\pa{\sup_{s \leq t}\int_0^sI\krull{x \leq X_u \leq y}dM_u}^p}\\
&\leq & K\EpX{}{\pa{\int_0^tI\krull{x \leq X_u \leq y}d\qvar{M}_u}^{p/2}}.
\eeqarno
where we use $K$ to denote any constant.
In order to estimate the last integral we now set $h=x-y$ and choose a
$C^2$-function $f$ with $f'' \geq 2I_{\left(x,y \right]}$ and $\abs{f'} \leq 2h$. Applying the Ito formula and
recalling that $\qvar{X} = \qvar{M}$ we obtain
\beqarno
\int_0^t I_{\left(x,y \right]}(X_s)d\qvar{M_s} &\leq & \half \int_0^t f''(X_s)d\qvar{M}_s =
f(X_t)-f(X_0) - \int_0^tf'(X_s)dX_s \\
&=&f(X_t)-f(X_0) - \int_0^tf'(X_s)dA_s - \int_0^tf'(X_s)dM_s\\
&\leq & 4Ch + \abs{\int_0^tf'(X_s)dM_s}.
\eeqarno
From the BDG inequality we now have
$$
\EpX{}{\pa{\sup_{s \leq t}\int_0^sf'(X_u)dM_u}^{p/2}}\leq \EpX{}{\pa{\int_0^t\abs{f'(X_s)}^2d\qvar{M}_s}^{p/4}}
\leq \pa{2Ch}^{p/2}.
$$
Choosing any $p > 2$ and using the Kolmogorov criterion from Theorem \ref{kc} gives
us the continuity result. \endproof
%%%%%%%%%%%%%%%%%%%%%%%%%%%%%%%%%%%%%%
%%%%%%%%%%%%%%%%%%%%%%%%%%%%%%%%%%%%%%%%%%%%%%%%%%%%%%%%%%%%%%%%%%%%%%%%%
\subsection{The Tanaka formulas} \label{tfb}
%%%%%%%%%%%%%%%%%%%%%%%%%%%%%%%%%%%%%%%%%%%%%%%%%%%%%%%%%%%%%%%%%%%%%%%
Given our definition of local time, the Tanaka formulas follow immediately.
%%%%%%%%%%%%%%%%%%%%
\bprop[The Tanaka formulas]
For any continuous semimartingale
$X$ we have the following relations.
\beqarno
\abs{X_t-x}&=&\abs{X_0-x}+\int_0^t \mbox{\em sgn}(X_s-x)dX_s + L_t^x(X),\\
\pa{X_t-x}^+&=&\pa{X_0-x}^+ +\int_0^t I\krull{X_s> x}dX_s + \half L_t^x(X),\\
\pa{X_t-x}^-&=&\pa{X_0-x}^- -\int_0^t I\krull{X_s \leq  x}dX_s + \half L_t^x(X).\\
\eeqarno
\eprop

\proof
The first formula is nothing more than Definition \ref{fp1}. The second
formula follows immediately from applying the first formula to the relation
$\abs{x}+x=2x^{+}$. The third formula follows in the same fashion.
%%%%%%%%%%%%%%%%%%%%%%%%%%%%%%%%%%%%%%%%%%%%%%%%%%%%%%%%%%%%%%%%%%%%%%%%%
\subsection{Local time as occupation time density} \label{odb}
%%%%%%%%%%%%%%%%%%%%%%%%%%%%%%%%%%%%%%%%%%%%%%%%%%%%%%%%%%%%%%%%%%%%%%%
Let $\map{f}{R}{R}$ be a convex function. A well known result then says that the
left hand derivative $f_-'(x)$ exists for every $x$. One can furthermore show that
$f_-'$ is left continuous and nondecreasing, so we can define a Borel measure $\mu_f$
on the real line by the prescription
$$
\mu_f\pa{\left[x,y \right)}=f_-'(y)-f_-'(x).
$$
We will also denote this measure by $f''(dx)$, and when $f \in C^2$ we will have
$\mu_f(dx)=f''(dx)=f''(x)dx$ where the last occurrence of $f''$ denotes the second order derivative. 
We now have the following generalization of the
Ito and Tanaka formulas.
%%%%%%%%%%%%%%%%%%%%%%%%%%%%%%
\bth[Extended Ito-Tanaka formula] \label{odbt1}
%%%%%%%%%%%%%%%%%%%%%%%%%%
Assume that $f$ is the difference between two convex functions, and that $X$ is a continuous semimartingale with local
time $L_t^x$. Then we have
\beq \label{odb1}
f(X_t)=f(X_0)+ \int_0^tf_-'(X)_sdX_S + \half \int_RL_t^xf''(dx)
\eeq
\eth

\proof 
When $f(x)=\abs{x}$ this is just the Tanaka formula with the measure $f''(dx)$  being a point
mass of size $2$ at $x=0$. If $f''$  consists of a finite number
of point masses, corresponding to an $f$ with a finite number of discontinuities
in $f'$, the formula follows easily from the Tanaka formulas and linearity. The
general case can then be proved by approximating a general $f''$ with a sequence
of measures with a finite number of point masses. See \cite{kall} for details. \endproof

We now have the following consequence of the extended Ito-Tanaka formula.
%%%%%%%%%%%%%%%%%%%%%%%%%%%%%%
\bth[Occupation time density]  \label{otd}
%%%%%%%%%%%%%%%%%%%%%%%%%%
For every nonnegative Borel function $\map{f}{R}{R_+}$ we have
\beq \label{odb2}
\int_0^tf(X_s)d\qvar{X}s=\int_Rf(x)L_t^xdx.
\eeq
\eth

\proof
Suppose that $f \in C$. Then we may interpret $f$ as $f=F''$ for a convex
function $F \in C^2$. Applying formula \cref{odb1} to $F(X_t)$  and comparing the result to
the standard Ito formula, shows that \cref{odb2} holds in this case. The general case
of a nonnegative Borel function $f$ then follows by a standard monotone class
argument. \endproof
%%%%%%%%%%%%%%%%%%%%%%%%%%%%%%%%%%%%%%%%%%%%%%%%%%%%%%%%%%%%%%%%%%%%%%%%%
\section{Notes to the literature} \label{n}
%%%%%%%%%%%%%%%%%%%%%%%%%%%%%%%%%%%%%%%%%%%%%%%%%%%%%%%%%%%%%%%%%%%%%%%
Full proofs of the main results (and much more) can be found in many textbooks
on stochastic calculus, such as \cite{cw},   \cite{kall}, \cite{kar-shr96}, \cite{P},  \cite{rev-yor}, and \cite{rw} . An extremely readable old classic
is \cite{mckean} which also contains a very nice discussion on regulated diffusions. Apart from a treatment
of the modern approach, \cite{kar-shr96} also contains an exposition of local time based on
the original Levy theory of excursions. An approach based on random walk
arguments is presented in \cite{mor-per}. The handbook \cite{bor-sal} contains a wealth of results.
%%%%%%%%%%%%%%%%%%%%%%%%%%%%%
\appendix
%%%%%%%%%%%%%%%%%%%%%%%%%%%%%%
%%%%%%%%%%%%%%%%%%%%%%%%%%%%%
\section{Stochastic Calculus for Continuous Martingales} \label{as}
%%%%%%%%%%%%%%%%%%%%%%%%%%%%%%
In this appendix we provide some basic concepts and facts concerning stochastic
integrals and the Ito formula for continuous semimartingales. We consider a
given filtered probability space $\filtspace$
 where $\op{F} = \krull{\Ft{t}}_{t \geq 0}$. Now consider
a process $X$ and a process property $E$. The property $E$ could be anything, like being a
martingale, being square integrable, being bounded, having finite variation etc. 
We now define the local version of $E$.
%%%%%%%%%%%%%%%%%%%%%%%%%%%%
\bdf
We say that $X$ has the property $E$ {\bf locally} if there exists an
increasing sequence of stopping times $\krull{\tau_n}_{n=1}^{\infty}$
with $\lim_{n \rightarrow \infty} \tau_n= \infty$ such that
the stopped process $X^{\tau_n}$ defined by
$$
X^{\tau_n}_t=X_{t \wedge \tau_n}
$$
has the property $E$ for each $n$. The sign $\wedge$ denotes the minimum operation.
\edf
%%%%%%%%%%%%%%%%%%%%%%%%%%

%%%%%%%%%%%%%%%%%%%%%%%%%%%%%
\subsection{Quadratic variation} \label{qv}
%%%%%%%%%%%%%%%%%%%%%%%%%%%%%%
Suppose that $X$ is a locally square integrable martingale with continuous trajectories.
Then $X^2$ will be a submartingale, and from the Doob-Meyer Theorem
we know that there exists a unique adaptive increasing process $A$ with $A_0=0$
such that the process
$$
X_t^2-A_t
$$
is a martingale, and one can show that $A$ is in fact continuous. The process $A$
above is very important for stochastic integration theory so we give it a name.
%%%%%%%%%%%%%%%%%%%%%%%%
\bdf
For any locally square integrable martingale $X$, the quadratic
variation process $\qvar{X}$ is defined by
$$
\qvar{X}_t=A_t,
$$
where $A$ is defined above.
\edf
If $X$ is a Wiener process $W$, then it is easy to see that $W_t-t$ is a martingale,
so we have the following trivial result.
%%%%%%%%%%%%%%%%%%%%%%%%%%%%%%%%%%%%%%
\begin{lemma}
For a Wiener process $W$ we have $qvar{W}_t=t$.
\end{lemma}
%%%%%%%%%%%%%%%%%%%%%%%%%%%%%%%%%%%%%%%%
The name ´´quadratic variation'' is motivated by the following result.
%%%%%%%%%%%%%%%%%%%%%%%%%%%%%%%%%%%%%%%%%%%%%%%
\bprop
Let $X$ be as above, consider a fixed $t >0$ and let $\krull{p_n}_{n=1}^{\infty}$
n=1
be a sequence of partitions of the interval $\bracket{0,t}$ such that the mesh (the longest
subinterval) of $p_n$ tends to zero as $n \rightarrow \infty$. Define $S_t^n$
$$
S_t^n=\sum_i\bracket{X\pa{t_{i+1}^n}-X\pa{t_{i}^n}}
$$
Then $S^n$ converges in probability to $\qvar{X}_t$.
\eprop
%%%%%%%%%%%%%%%%%%%%%%%%%%%%%%%%%%
This results motivates the formal expression
\beq \label{qv1}
d\qvar{X}_t=\pa{dX_t}^2.
\eeq
and we note that for the case when $X$ is a Wiener process $W$, this is the usual ``multiplication rule'' $\pa{dWt}2 = dt$.

We end this section by quoting the important Burkholder-Davis-Gundy (BDG) inequality.
This inequality shows how the maximum of a local martingale is
controlled by its quadratic variation process.
%%%%%%%%%%%%%%%%%%%%%%%%%%%%%%
\bth[Burkholder-Davis-Gundy] \label{bdg}
%%%%%%%%%%%%%%%%%%%%%%%%%%%%%%%%%%%%%%
For every $p >0$ there exists a
constant $C_p$ such that
\beq \label{qv2}
\EpX{}{\pa{M_t^{\star}}^p}\leq C_p\EpX{}{\qvar{M}_t^{p/2}}
\eeq
for every continuous local martingale $M$, where $M^{\star}$ defined by
$$
M_t^{\star}=\sup_{s \leq t}\abs{M_s}.
$$
\eth
%%%%%%%%%%%%%%%%%%%%%%%%%%%%%%%%%%%%%
\subsection{Stochastic integrals}\label{si}
%%%%%%%%%%%%%%%%%%%%%%%%%%%%%%%%%%%%%%%%%%%%%
Assume that $M$ is a continuous square integrable martingale. Then one can quite easily
define the stochastic integral 
$$
\int_0^tg_sdM_s
$$
almost exactly along the lines of the construction of the usual Ito integral. The
basic properties of the integral are as follows.
%%%%%%%%%%%%%%%%%%%%%%%%%%%%%
\bprop
%%%%%%%%%%%%%%%%%%%%%%%%%%
 Let M be a square integrable martingale, and let g be an
adapted process satisfying the integrability condition
$$
\EpX{}{\int_0^tg_s^2d\qvar{M}_s} < \infty
$$
for all $t$. Then the integral process $g \star M$ defined by
$$
\pa{g \star M}_t=\int_0^tg_sdM_s,
$$
is well defined and has the following properties.
\bei
\item
The process $g \star M$ is a square integrable martingale with continuous trajectories.
\item
$$
\EpX{}{\pa{\int_0^tg_sdM_s}^2}=\EpX{}{\int_0^tg_s^2d\qvar{M}_s}
$$
\item 
We have
\beq \label{si5}
d\qvar{g \star M}_t=g_t^2 d\qvar{M}_t.
\eeq
\ei
\eprop
The proof of the first two items above are very similar to the standard Ito
case. The formula \cref{si5} follows from \cref{qv1} and the fact that $d(g \star M)_t=g_tdM_t$.
The stochastic integral above can quite easily be extended to a larger class
of processes. Note that a continuous local martingale is also locally square
integrable.
%%%%%%%%%%%%%%%%%%
\bdf 
For a continuous local martingale $M$, we define $L^2(M)$ as the
class of adapted processes $g$ such that
$$
\int_0^tg_s^2d\qvar{M}_s < \quad P-a.s.
$$
for all $t \geq 0$.
\edf
We now have the following result.
%%%%%%%%%%%%%%%%%%%
\bprop
 Let $M$ be a local martingale, and let $g$ be an adapted process
in $L^2(M)$. Then the integral process $g \star M$ defined by
$$
\pa{g \star M}_t=\int_0^tg_sdM_s,
$$
is well defined and has the following properties.
\bei
\item
The process $g \star M$ is a local martingale with continuous trajectories.
\item
We have
\beq \label{si10}
d\qvar{g \star M}_t=g_t^2 d\qvar{M}_t.
\eeq
\ei
\eprop
%%%%%%%%%%%%%%%%%%%%%%%%%%%%%%%%%%%%%
\subsection{Semimartingales and the Ito formula} \label{sit}
%%%%%%%%%%%%%%%%%%%%%%%%%%%%%%%%%%%%%%%%%%%%%
We start by defining the semimartingale concept.
\bdf
If a process $X$ has the form
$$
X_t=A_t+ M_t
$$
where $A$ is adapted with $A_0=0$, continuous and of (locally) bounded variation,
and $M$ is a local martingale, then we say that $X$ is a {\bf semimartingale}. If $A$
can be written on the form
$$
dA_t=\mu_tdt,
$$
then we say that $X$ is a special semimartingale.
\edf
%%%%%%%%%%%%
We have the following uniqueness result.
%%%%%%%%%%%
\begin{lemma}
The decomposition $X = M + A$ is unique.
\end{lemma}
Given the stochastic integral defined earlier, we see that semimartingales
are natural integrators.
\bdf
Let $X$ be a semimartingale as above. The class $L(X)$ is defined
as
$$
L(X) = L(M) \cap L^1(\abs{dA}).
$$
where $L^1(\abs{dA})$ is the class of adapted processes $g$ such that
$$
\int_0^t\abs{g_s}\abs{dA}_s < \infty
$$
for all $t$.
For $ g \in L(X)$ we define the stochastic integral by
$$
\int_0^tg_sdX_s=\int_0^tg_sdA_s+ \int_0^tg_sdM_s.
$$
The quadratic variation of $X$ is defined by
$$
d\qvar{X}_t=d\qvar{M}_t.
$$
\edf
It is easy to show that any function of bounded variation has zero quadratic
variation, so \cref{qv1} can be extended to
$$
d\qvar{X}_t = \pa{dX_t}^2.
$$
Very much along the lines of the informal rules $(dt)^2=0$, and $dt \cdot dW_t=0$ in
standard Ito calculus we can thus motivate the informal multiplication table
\beqarno
\pa{dX_t}^2&=&d\qvar{X}_t,\\
\pa{dA_t}^2&=&0,\\
dA_t\cdot dM_t&=&0.
\eeqarno
Let us now consider a semimartingale $X$ as above, and assume that $\map{F}{R}{R}$
is a smooth function. Inspired by the standard Ito formula we would then
intuitively write
$$
dF(X_t)=F'(X_t)dX_T+ \half F''(X_t)\pa{X_t}^2.
$$
Given the multiplication table above,
we also have
$$
\pa{dX_t}^2=\pa{dA_t+ dM_t}^2=\pa{dM_t}^2,
$$ 
so from \cref{qv1} we would finally expect to obtain
$$
\pa{dX_t}^2=d\qvar{M}_t.
$$
This intuitive argument can in fact be made precise and we have the following
extension of the Ito formula.
%%%%%%%%%%%%%%%%%%%%%%%%
\bth[The Ito Formula]
Assume that $X$ is a semimartingale and
that $F(t,x)$ is $C^{1,2}$. We then have the Ito formula
$$
dF(t,X_t)=\dfdx{F}{t}(t,X_t)dt + \dfdx{F}{x}(t,X_t)dX_t + \half \ddfdxx{F}{x}(t,X_t)d\qvar{X_t}.
$$
\eth

We end the section by quoting the extremely useful dominated stochastic convergence
theorem.
%%%%%%%%%%%%%%
\bth[Stochastic dominated convergence] \label{dc}
Let X be a continuous semimartingale
and let $Z,Y,Y^1,Y^2,\ldots$  be processes in $L(X)$. Assume that the following
conditions hold for all $t$.
\beqarno
Y_t^n & \rightarrow & Y_t, \quad P-a.s.\\
\abs{Y_t^n} &\leq & \abs{Z_t},\quad \mbox{for all $n$.}
\eeqarno
Then we have
$$
\sup_{s \leq t}\abs{\int_0^sY_u^n dX_u-\int_0^sY_u dX_u}\stackrel{P}{\rightarrow} 0.
$$
\eth
%%%%%%%%%%%%%%%%%%%%%%%%%%%%%%%%%%%%%
\subsection{The Levy characterization of Brownian motion} \label{lev}
%%%%%%%%%%%%%%%%%%%%%%%%%%%%%%%%%%%%%%%%%%%%% 
Let $W$ be a Wiener process. It is very easy to see that $W$ has the following
properties
\bei
\item
$W$ is a continuous martingale.
\item
The process$W_t^2-t$ is a martingale.
\ei
Surprisingly enough, these properties completely characterize the Wiener process
so we have the following result.
%%%%%%%%%%%%%%%%%%%%%%%
\bprop[Levy] 
Assume that the process $X$ has the following properties.
\bei
\item
$X$ has continuous trajectories and $X_0=0$.
\item
$X$ is a martingale.
\item
The process$X_t^2-t$ is a martingale.
\ei
Then $X$ is a standard Wiener process.
\eprop
We have the following easy corollary.
%%%%%%%%%%%%%%%%
\begin{corollary}
Assume that $X$ is a continuous martingale with quadratic variation
given by
$$
d\qvar{X}_t=dt.
$$
Then $X$ is a standard Wiener process.
\end{corollary}

\proof
 We apply the Ito formula to $X^2$ to obtain
$$
d\pa{X_t^2}=2X_tdX_t + d\qvar{X}_t=2X_tdX_t + dt
$$
and, since $X$ is a martingale, the term $X_tdX_t$ is a martingale increment.
\endproof

We also have the following easy and useful consequence of the Levy characterization.
%%%%%%%%%%%%%%
\bprop \label{lev10}
%%%%%%%%%%%%%%%%%%%%
Assume that $X$ is a continuous martingale such that
$$
d\qvar{X}_t=\sigma_t^2dt
$$
for some process$\sigma >0$. Then the process $W$, defined by
$$
W_t=\frac{1}{\sigma_t}dX_t
$$
is a standard Wiener process, so we can write
$$
dX_t=\sigma_tdW_t
$$
where $W$ is standard Wiener.
\eprop

\proof
Since $W$ is a stochastic integral with respect to a martingale we conclude
that $W$ is a martingale. Using \cref{si10} we obtain
$$
d\qvar{W}_t=\frac{1}{\sigma_t^2}d\qvar{X}_t=dt. \quad \endproof
$$

%%%%%%%%%%%%%%%%%%%%%%%%%%%%%%%%%%%%%
\section{A random time change} \label{rt}
%%%%%%%%%%%%%%%%%%%%%%%%%%%%%%%%%%%%%%%%%%%%% 
In the section we will investigate how the structure of a semimartingale changes
when we perform a random time change.
\bdf
A random time change is a process $\krull{C_t: \  t \geq 0}$ such that
\bei
\item
$C_0=0$
\item
$C$ has strictly increasing continuous trajectories.
\item
For every fixed $t \geq 0$, the random time $C_t$ is a stopping time.
\ei
We say that $C$ is smooth if the trajectories are differentiable.
\edf
%%%%%%%%%%%%%%%%%%%%%%%
In a more general theory of random time changes, we only require that $C$ is
non decreasing, but the assumption of strictly increasing continuous trajectories
makes life much easier for us.
\bdf
If $X$ is an adapted process and $C$ is a random time change,
the process $X^C$ is defined as
$$
X_t^C=X_{C_t}.
$$
\edf
We remark that if $X$ is continuous and adapted to the filtration $\op{F}=\bigfilt$
then $X^C$ is adapted to the filtration $\op{F}^C=\krull{\Ft{C_t}}_{t \geq o}$.
The following extremely useful result tells us how the quadratic variation
changes after a time change.
%%%%%%%%%%%%%%%%%
\bprop \label{rt5}
%%%%%%%%%%%%%%%%%%%%%
Let $X$ be a local $\op{F}$-martingale and let $C$ be a continuous
random time change. Then the following hold.
\bei
\item
$X^C$ is a local $\op{F}^C$-martingale.
\item
The quadratic variation of $X^C$ is given by
$$
\qvar{X^C}_t=\qvar{X}_t^C,
$$
or equivalently
$$
\qvar{X^C}_t=\qvar{X}_{C_t}.
$$
\ei
\eprop
%%%%%%%%%%%%%
\proof
The first part follows from the optional sampling theorem. For the
second part we assume WLOG that $X$ is a square integrable martingale and
recall that $A_t=\qvar{X}_t$ is the unique continuous process of bounded variation
such that
$$
X_t^2-A_t
$$
is an $\op{F}$-martingale. From the optional sampling theorem we conclude that $Z$, defined by
$$
Z_t=X_{C_t}^2-A_{C_t}
$$
is an $\op{F}^C$-martingale. On the other hand we trivially have
$$
A_{C_t}=A_t^C, \quad \mbox{and} \quad \pa{X^2}_{C_t}=\pa{X^C}_t^2,
$$
so the process
$$
\pa{X^C}_t^2-A_t^C
$$
is an $\op{F}^C$-martingale. 
\endproof

We will now investigate the effect of a random time change to a stochastic
integral. Let us consider a local martingale $X$ of the form
$$
X_t=\int_0^t\sigma_sdW_s,
$$
and a random time change $C$. We know that
$$
\qvar{X}_t=\int_0^t\sigma_s^2ds
$$
so from Proposition \ref{rt5} we see that 
\beq \label{rt15}
\qvar{X^C}_t=\int_0^{C_t}\sigma_s^2ds,
\eeq
so if the time change is smooth we have
\beq \label{rt20}
\qvar{X^C}_t=\sigma_{C_t}C_t'dt.
\eeq
From this we have the following important result which follows directly from \cref{rt20} and Proposition \ref{lev10}.
%%%%%%%%%%%%%%%
\bprop \label{rt21}
%%%%%%%%%%%%%%%%%%%%%%%%%%%%%
Let $X$ be defined by
$$
X_t=\int_0^t\sigma_sdW_s
$$
and consider a smooth time change $C$. Then we can write
$$
dX_t^C=\sigma_{C_t}\sqrt{C_t'}d\bar{W}_t
$$
where $\bar{W}$ is standard Wiener.
\eprop
%%%%%%%%%%%%%%%%
The most commonly used special cases are the following.
%%%%%%%%%%%%%%%%%
\bprop \label{rt22}
%%%%%%%%%%%%%%%%%%%%%%
With $X$ as above, define the random time change $C$ by the
implicit relation
$$
\int_0^{C_t}\sigma_s^2ds=t.
$$
i.e.
$$
C_t=\inf \krull{s \geq 0: \ \int_0^s\sigma_u^2du=s}
$$
Then the process $X_t^C=X_{C_t}$ is a standard Wiener process.
\eprop
%%%%%%%%%%
\proof
Direct differentiation shows that
$$
C_t'=\frac{1}{\sigma_{C_t}^2}
$$
and we can now use Proposition \ref{rt21}. \endproof
%%%%%%%%%%%%%%%%%%%%%%%%%%%%%%%%%%%
\bprop \label{rt23}
%%%%%%%%%%%%%%%%%%%%%%%%%%
Assume that ˜$\bar{W}$ is standard Wiener, and assume that $\sigma > 0$. Define the random time
change $C$ by the implicit relation
$$
\int_0^{C_t}\sigma_s^{-2}ds
$$
i.e.
$$
C_t = \inf \krull{s \geq 0:\ \int_0^s \sigma_u^{-2}du=t}.
$$
Then the process $\bar{W}^C$ can be written as
˜$$
\bar{W}_t^C=\int_0^t \sigma_s dW_s,
$$
where $W$ is standard Wiener.
\eprop
%%%%%%%%%%%%%%%%%%%%%%

We may in fact push this analysis a bit further. Let us therefore consider
an SDE of the form
$$
dX_t=\mu(X_t)dt + \sigma (X_t)dW_t,
$$
and let us define the time change $C$ by
$$
\int_0^{C_t}g^2(X_s)ds=t
$$
for some deterministic function $g >0$. This transformation is in fact smooth
and we have
$$
C_t'=\frac{1}{g^2\pa{X_t^C}}.
$$
We  have
$$
X_t^C=X_0+ \int_0^t\mu(X_s)ds + \int_0^t\sigma (X_s)dW_s
$$
so we obtain
$$
X_t=X_0+ \int_0^{C_t}\mu(X_s)ds + \int_0^{C_t}\sigma (X_s)dW_s
$$
Defining $Y$ by
$$
Y_t=\int_0^{C_t}\mu (X_s)ds
$$
we obtain
$$
dY_t=\mu \pa{X_t^C}C_t'dt = \frac{\mu \pa{X_t^C}}{g^2\pa{X_t^C}}dt.
$$
Defining $Z$ by
$$
Z_t=\int_0^{C_t}\sigma (X_s)dW_s
$$
and using Proposition \ref{rt5} we obtain
$$
\qvar{Z}_t=\int_0^{C_t}\sigma^2(X_s)ds
$$
so on differential form we have
$$
d\qvar{Z}_t=\sigma^2\pa{X_t^C}C_t'dt=\frac{\sigma^2\pa{X_t^C}}{g^2 \pa{X_t^C}}
$$
It now follows from Proposition \ref{lev10} that we can write
$$
dZ_t=\frac{\sigma\pa{X_t^C}}{g \pa{X_t^C}}d\bar{W}_t
$$
where $\bar{W}$ is standard Wiener. We have thus proved the following result.
%%%%%%%%%%%%%%%%%%%%%%%%%%%
\bprop \label{rt25}
%%%%%%%%%%%%%%%%%%%%%%%%%%%
Consider the SDE
$$
dX_t=\mu(X_t)dt + \sigma (X_t)dW_t,
$$
Consider furthermore a function $g >0$ and and a random time change of the
form 
$$
\int_0^{C_t}g^2(X_s)ds=t
$$
or equivalently
$$
C_t=\inf \krull{s \geq 0:\ \int_0^{C_t}g^2 \pa{X_s}ds=t}.
$$
Defining $Y$ by $Y=X^C$, the process $Y$ will satisfy the SDE
$$
dY_t=\mu_Y(Y_t)dt + \sigma_Y(Y_t)d\bar{W}_t
$$
where
$$
\mu_Y(y)=\frac{\mu(y)}{g^2(y)},\quad \sigma_Y(y)=\frac{\sigma(y)}{g(y)}
$$
and $\bar{W}$ is standard Wiener.
\eprop
%%%%%%%%%%%%%%%%%%%%%%%%%%%%%%%%%%%%%%%%%%%%%%%
\section{The Kolmogorov continuity criterion} \label{kc}
%%%%%%%%%%%%%%%%%%%%%%%%%%%%%%%%%%%%%%%%%%%%%%%
Whenever you want to prove continuity for some (possibly multi-indexed) process,
the Kolmogorov criterion is often the first choice.
%%%%%%%%%%%%%%%%%%%%
\bth[Kolmogorov continuity criterion]
 Let $\map{X}{\Omega \times R^D}{S}$ 
be a process indexed by $R^d$ and taking values in some complete metric space $(S, \rho)$
 Assume that for some positive real numbers $a$, $b$,and $C$ we have
$$
\EpX{}{\krull{\rho \pa{X_t,X_s}}^a}\leq C\norm{s-t}{}^{d+b},\quad \mbox{for all $s,t \in R^D$}
$$
Then $X$ has a continuous version.
\eth
%%%%%%%%%%%%%%%%%%%%%%%%%
The two most common choices of $S$ are $S=R^n$ or $S=C(R)$ with the uniform
topology.
%%%%%%%%%%%%%%%%%%%%%%%%%%%%%%% 
\bibliography{REF}
%%%%%%%%%%%%%%%%%%%%%%%%%%%%%
%%%%%%%%%%%%%%%%%%%%%%%%%%%%%%%%%%%
\end{document}